\documentclass[12pt,reqno]{amsart}
\usepackage{fullpage}
\usepackage{times}
\usepackage{multirow,graphicx}

\renewcommand{\labelenumi}{(\alph{enumi})}

\begin{document}

\title{An annotated bibliography of work related to gender in science}
\author[Greg Martin]{Greg Martin} 
\maketitle
\thispagestyle{empty}

The purpose of this manuscript is to gather together a large amount of source material pertaining to women in mathematics, from studies of girls in elementary school through data on females winning prizes for mathematical research. Along the way, we have also gathered a large amount of material from the psychology and sociology literature on implicit biases more generally, particularly pertaining to gender. This source material was then used to support the writing of the article~\cite{Martin}. We have tried to refer to primary research literature whenever possible, although we have also included well-written blog posts, organizational web sites, self-published articles by research organizations, and even a YouTube video.

Each bibliography entry is accompanied by some remarks summarizing its content and representative quotes from the articles themselves. We have followed standard practice when including these quotes, with the following exception: where the original quote has included citations to other work, or supporting statistics such as $p$-values, we have omitted these annotations to enhance the clarity of the quote. Nevertheless, much of the work in this bibliography contains a large number of further references to the relevant research literature.

The annotated bibliography is thus reasonable for browsing; but for those looking to find source material for particular aspects of this issue, we hope the following categories will be of some use:

\medskip
{
\smaller\smaller
\begin{itemize}
\item Hypotheses on the causes of underrepresentation of women in science: \cite{Bar}, \cite{EHL}, \cite{GMSZ}, \cite{KM}, \cite{LHPL}, \cite{Tyson}
\item Declining gender gap in math achievement by girls and boys: \cite{AGKM}, \cite{EHL}, \cite{GMSZ}, \cite{HLLEW}, \cite{HM}, \cite{KM}, \cite{OECD}
\item Overemphasis of, and problems with, standardized tests: \cite{AAUW}, \cite{EHL}, \cite{HLLEW}, \cite{NV}
\item Children with extremely high math achievement: \cite{AGKM}, \cite{GMSZ}, \cite{HLLEW}, \cite{HM}, \cite{NV}
\item Role of culture's gender inequity in math achievement by girls and boys: \cite{EHL}, \cite{GMSZ}, \cite{HM}, \cite{KM}
\item Biases in primary school classrooms: \cite{AAUW}, \cite{BGRL}, \cite{LHPL}, \cite{NV}, \cite{W}
\item Intelligence as fixed trait vs.\ intelligence as malleable quality: \cite{GRD}, \cite{MD}, \cite{U}
\item Explicit sexism: \cite{Bar}, \cite{GF:C}, \cite{Hess}, \cite{Tyson}, \cite{WISELI:o}
\item Implicit biases: \cite{Bac}, \cite{BMN}, \cite{BHKM}, \cite{DK}, \cite{DPAHV}, \cite{G}, \cite{GR}, \cite{GBN}, \cite{HWFT}, \cite{MAC}, \cite{MDBGH}, \cite{RSZ}, \cite{RP}, \cite{SAR}, \cite{TP}, \cite{UC}, \cite{WISELI:o}, \cite{WISELI:S}
\item Incognizance of our biases, and the illusion of meritocracy: \cite{Bac}, \cite{G}, \cite{GBN}, \cite{KGH}, \cite{PGR}, \cite{Ries}, \cite{UC}, \cite{VS}
\item Broader societal gender-based problems: \cite{AAUW}, \cite{AGKM}, \cite{Bac}, \cite{Bar}, \cite{BMN}, \cite{CI}, \cite{EK}, \cite{HWHH}, \cite{HWFT}, \cite{Hess}, \cite{KGH}, \cite{LHPL}, \cite{RC}
\item Gender-based personality expectations: \cite{BLGS}, \cite{BBL}, \cite{CI}, \cite{DK}, \cite{EK}, \cite{EIK}, \cite{HWFT}, \cite{N}, \cite{R}, \cite{Snyder}
\item Gender-based differential in self-concept: \cite{NV}, \cite{RSZ}, \cite{SK}, \cite{UC}
\item Effect of parenting, childcare, flexible schedules on one's career: \cite{BGS}, \cite{CBP}, \cite{DK}, \cite{FP:C}, \cite{Lalin}, \cite{Munsch}, \cite{RP}, \cite{Sa}
\item Impostor phenomenon: \cite{CI}, \cite{Kap}
\item Stereotype threat: \cite{FP:G}, \cite{GRD}, \cite{LHPL}, \cite{NV}, \cite{SSQ}
\item Different standards for women and men regarding leadership, persuasion, and negotiation: \cite{BLGS}, \cite{BBL}, \cite{BHKM}, \cite{CE}, \cite{DK}, \cite{EK}, \cite{EIK}, \cite{HSR}, \cite{R}
\item Female speakers at conferences: \cite{AMSsurvey}, \cite{AWM}, \cite{CH}, \cite{DK}, \cite{FP:C}, \cite{FP:G}, \cite{F}, \cite{GF:C}, \cite{G}, \cite{H}, \cite{IYH}, \cite{JSConf}, \cite{OM}, \cite{P}, \cite{Ries}, \cite{RM}, \cite{St}, \cite{Su}, \cite{VS}
\item Gender-based differential in teaching evaluations: \cite{Kas}, \cite{N}, \cite{SK}
\item Gender-based differential in award winners, grants, promotion/tenure: \cite{CGFSH}, \cite{CMR}, \cite{HWFT}, \cite{KGH}, \cite{SS}, \cite{SAR}, \cite{WISELI:o}
\item Gender-based biases in evaluation and selection processes: \cite{BGS}, \cite{CGFSH}, \cite{DE}, \cite{GR}, \cite{H}, \cite{HSR}, \cite{IYH}, \cite{Kas}, \cite{KGH}, \cite{MAC}, \cite{MDBGH}, \cite{Munsch}, \cite{RSZ}, \cite{SS}, \cite{SK}, \cite{Snyder}, \cite{TP}, \cite{VMC}, \cite{WW}, \cite{WISELI:S}
\item Action items and recommendations for addressing underrepresentation of women: \cite{AAUW}, \cite{AGKM}, \cite{CGFSH}, \cite{CH}, \cite{DK}, \cite{FP:C}, \cite{FP:G}, \cite{GF:C}, \cite{GF:T}, \cite{IMA}, \cite{IYH}, \cite{JSConf}, \cite{Lalin}, \cite{NSF} \cite{OM}, \cite{RM}, \cite{St}, \cite{VS}, \cite{W}, \cite{WISELI:o}, \cite{WISELI:S}
\item Further references: \cite{AAUW}, \cite{AGKM}, \cite{DK}, \cite{DE}, \cite{EHL}, \cite{EIK}, \cite{FP:G}, \cite{GRD}, \cite{GMSZ}, \cite{KM}, \cite{KGH}, \cite{N}, \cite{PGR}, \cite{R}, \cite{W}, \cite{WISELI:o}
\end{itemize}
} 
\medskip

This manuscript also includes an appendix containing tables of data from the 2014 ICM and the 2014 Joint Meetings of the AMS and MAA, listing the various sessions and the numbers of female speakers and total speakers (and the same for organizers, where given). Some data from mathematics prizes will also be included.

As we write in \cite{Martin}: ``We have made the conscious choice to include only initials and last names in the bibliography and in both manuscripts. We have observed a tendency to be curious about the gender of the authors of the research referred to herein, and perhaps to involuntarily wonder how the authors' gender should affect our evaluation of their conclusions. These reflexive speculations, we believe, tellingly illuminate the depth to which these implicit biases about gender are ingrained in us, even though we rationally know that possessing one gender or another does not affect a person's objectivity. Being socialized to have biases is not our fault; but preventing our biases from negatively affecting the world around us is nonetheless our responsibility.''

\section*{Acknowledgments}

We thank W. Miao for gathering the data appearing in the appendix, as well as for locating copies of several of the papers in this bibliography; we also thank J. Bryan for performing the statistical analysis described in the appendix.

\section*{Appendix: Gender data for conferences, prizes, and editorial boards}

This section contains all of the data used to calculate the statistics reported in Section 1.2 of~\cite{Martin}; we are grateful to W. Miao for gathering this data. In the tables below, the column heading {\bf FS} counts the number of female speakers (in which we include also panelists, prizewinners, and editors), {\bf TS} the total number of speakers, {\bf FO} the number of female organizers (in which we include also moderators, chairs, and introducers), and {\bf TO} the total number of organizers. Genders were determined via internet searches. Asterisks (*, **, ***) represent one, two, or three speakers whose genders could not be determined; these speakers were not counted in the total number of speakers.

{
\smaller\smaller

\bigskip\bigskip\centerline{\bf Joint AMS/MAA Winter Meeting, Baltimore, January 2014} \bigskip

\centerline{MAA Invited Paper Sessions}\bigskip

\centerline{
\begin{tabular}{|l||c|c||c|c|}\hline
\bf Session & \bf FS & \bf TS & \bf FO & \bf TO \\ \hline\hline
Continuing Influence of Paul Erd\H os in Number Theory & 1 & 6 & 0 & 2 \\ \hline
Graphs Don't Have to Lie Flat: The Shape of Topological Graph Theory & 2 & 4 & 1 & 2 \\ \hline
Mathematics and Effective Thinking & 5 & 11 & 0 & 1 \\ \hline
Six Crash Courses on Mapping Class Groups & 2 & 7 & 0 & 2 \\ \hline
Uniform Distribution, Discrepency, and Related Fields & 0 & 6 & 1 & 2 \\ \hline
Unreasonable Effectiveness of Modern Mathematics & 1 & 4 & 1 & 2 \\ \hline
\end{tabular}
}

\newpage

\centerline{Invited Addresses}\bigskip

\centerline{
\begin{tabular}{|l||c|c|}\hline
\bf Category & \bf FS & \bf TS \\ \hline\hline
Joint Invited Addresses & 0 & 3 \\ \hline
AMS Invited Addresses & 1 & 7 \\ \hline
MAA Invited Addresses & 3 & 7 \\ \hline
Invited Addresses of Other Organizations & 3 & 10 \\ \hline
\end{tabular}
}

\bigskip\bigskip

\centerline{AMS Special Sessions}\bigskip

\centerline{
\begin{tabular}{|l||c|c||c|c|}\hline
\bf Session & \bf FS & \bf TS & \bf FO & \bf TO \\ \hline\hline
Accelerated Advances in Higher Order Invexities/Univexities with Applications & \multirow2*{2} & \multirow2*{10} & \multirow2*{0} & \multirow2*{2} \\
\; to Optimization and Mathematical Programming &&&& \\ \hline
Advances in Analysis and PDEs	 & 5 & 14 & 0 & 2 \\ \hline
Algebraic Geometry* & 1 & 11 & 0 & 2 \\ \hline
Algebraic Structures Motivated by Knot Theory & 11 & 22 & 1 & 5 \\ \hline
Algebraic and Analytic Aspects of Integrable Systems and Painlev\'e Equations & 3 & 14 & 0 & 3 \\ \hline
Analytic Number Theory & 6 & 24 & 0 & 3 \\ \hline
Applied Harmonic Analysis: Large Data Sets, Signal Processing, and Inverse & \multirow2*{2} & \multirow2*{12} & \multirow2*{0} & \multirow2*{3} \\ 
\; Problems &&&& \\ \hline
Banach Spaces, Metric Embeddings, and Applications & 1 & 14 & 1 & 2 \\ \hline
Big Data: Mathematical and Statistical Modeling, Tools, Services, and Training & 2 & 6 & 0 & 1 \\ \hline
Categorical Topology & 0 & 10 & 0 & 2 \\ \hline
The Changing Education of Preservice Teachers in Light of the Common Core & 12 & 19 & 2 & 4 \\ \hline
Classification Problems in Operator Algebras & 2 & 17 & 0 & 2 \\ \hline
Communication of Mathematics via Interactive Activities & 4 & 10 & 0 & 2 \\ \hline
Complex Dynamics, I (a Mathematics Research Communities Session) & 9 & 20 & 1 & 3 \\ \hline
Computability in Geometry and Topology & 4 & 15 & 0 & 2 \\ \hline
De Bruijn Sequences and Their Generalizations & 2 & 12 & 0 & 2 \\ \hline
Deformation Spaces of Geometric Structures on Low-Dimensional Manifolds & 4 & 17 & 2 & 4 \\ \hline
Difference Equations and Applications & 2 & 14 & 0 & 1 \\ \hline
Dispersive and Geometric Partial Differential Equations & 4 & 18 & 0 & 3 \\ \hline
Ergodic Theory and Symbolic Dynamics & 6 & 22 & 1 & 2 \\ \hline
Fractal Geometry: Mathematics of Fractals and Related Topics & 2 & 20 & 0 & 4 \\ \hline
Fractional, Stochastic, and Hybrid Dynamic Systems with Applications & 4 & 20 & 1 & 3 \\ \hline
Geometric Applications of Algebraic Combinatorics & 10 & 24 & 2 & 2 \\ \hline
Geometric Group Theory, I (a Mathematics Research Communities Session) & 5 & 20 & 1 & 4 \\ \hline
Global Dynamics and Bifurcations of Difference Equations & 3 & 15 & 0 & 2 \\ \hline
Graph Theory: Structural and Extremal Problems & 1 & 22 & 0 & 2 \\ \hline
Heavy Tailed Probability Distributions and Their Applications & 0 & 12 & 0 & 2 \\ \hline
Highlighting Achievements and Contributions of Mathematicians of the African & \multirow2*5 & \multirow2*{14} & \multirow2*1 & \multirow2*2 \\
\; Diaspora &&&& \\ \hline
History of Mathematics & 8 & 22 & 2 & 3 \\ \hline
Homological and Characteristic p Methods in Commutative Algebra & 8 & 22 & 0 & 3 \\ \hline
Homotopy Theory & 7 & 22 & 0 & 5 \\ \hline
Hyperplane Arrangements and Applications* & 2 & 19 & 0 & 3 \\ \hline
Logic and Probability & 3 & 20 & 2 & 4 \\ \hline
Mathematics and Mathematics Education in Fiber Arts & 8 & 12 & 2 & 2 \\ \hline
Mathematics in Natural Resource Modeling & 6 & 15 & 2 & 2 \\ \hline
Mathematics of Computation: Differential Equations, Linear Algebra, and & \multirow2*5 & \multirow2*{15} & \multirow2*1 & \multirow2*2 \\
\; Applications* &&&& \\ \hline
\end{tabular}
}

\newpage

\centerline{AMS Special Sessions (continued)}\bigskip

\centerline{
\begin{tabular}{|l||c|c||c|c|}\hline
\bf Session & \bf FS & \bf TS & \bf FO & \bf TO \\ \hline\hline
My Favorite Graph Theory Conjectures & 3 & 19 & 1 & 2 \\ \hline
Nineteenth Century Algebra and Analysis & 2 & 10 & 1 & 3 \\ \hline
Nonlinear Systems: Polynomial Equations, Nonlinear PDEs, and Applications & 0 & 20 & 0 & 1 \\ \hline
Outreach for Mathematically Talented Youth & 9 & 15 & 2 & 3 \\ \hline
Progress in Free Analysis and Free Probability & 2 & 15 & 0 & 2 \\ \hline
Quantum Walks, Quantum Computation, and Related Topics*** & 0 & 14 & 0 & 4 \\ \hline
Random Matrices: Theory and Applications & 1 & 10 & 0 & 2 \\ \hline
Reaction Diffusion Equations and Applications & 4 & 16 & 0 & 2 \\ \hline
Recent Advances in Homogenization and Model Reduction Methods for Multiscale & \multirow2*3 & \multirow2*{16} & \multirow2*1 & \multirow2*2 \\
\; Phenomena &&&& \\ \hline
Recent Progress in Geometric and Complex Analysis & 3 & 19 & 0 & 3 \\ \hline
Recent Progress in Multivariable Operator Theory & 7 & 22 & 0 & 2 \\ \hline
Recent Progress in the Langlands Program & 0 & 18 & 0 & 2 \\ \hline
Regularity Problem for Nonlinear PDEs Modeling Fluids and Complex Fluids, I (a & \multirow2*5 & \multirow2*{20} & \multirow2*0 & \multirow2*4 \\
\; Mathematics Research Communities Session) &&&& \\ \hline
Representation Theory of p-adic Groups and Automorphic Forms & 1 & 8 & 0 & 2 \\ \hline
Research in Mathematics by Undergraduates and Students in Post-Baccalaureate & \multirow2*{16} & \multirow2*{30} & \multirow2*1 & \multirow2*6 \\
\; Programs &&&& \\ \hline 
Set-Valued Optimization and Variational Problems with Applications & 4 & 16 & 1 & 4 \\ \hline
Symplectic and Contact Structures on Manifolds with Special Holonomy & 2 & 14 & 1 & 3 \\ \hline
Topological Graph Theory: Structure and Symmetry & 1 & 22 & 0 & 2 \\ \hline
Trends in Graph Theory & 6 & 16 & 1 & 1 \\ \hline
Tropical and Nonarchimedean Analytic Geometry, I (a Mathematics Research & \multirow2*4 & \multirow2*{20} & \multirow2*1 & \multirow2*3 \\
\; Communities Session) &&&& \\ \hline
The Ubiquity of Dynamical Systems & 3 & 9 & 1 & 2 \\ \hline
\end{tabular}
}

\bigskip\bigskip
\centerline{MAA General Contributed Paper Sessions}
\centerline{(2 female organizers, 3 total organizers)}\bigskip

\centerline{
\begin{tabular}{|l||c|c|}\hline
\bf Session & \bf FS & \bf TS \\ \hline\hline
Assessment and Outreach & 3 & 6 \\ \hline
Assorted Topics I \& II & 6 & 23 \\ \hline
Calculus & 4 & 8 \\ \hline
History and Philosophy of Mathematics & 3 & 10 \\ \hline
Interdisciplinary Topics & 4 & 4 \\ \hline
Mathematics Education I, II, \& III & 24 & 47 \\ \hline
Mathematics and Technology & 4 & 9 \\ \hline
Modeling and Applications of Mathematics I, II, \& III* & 20 & 43 \\ \hline
Probability and Statistics I \& II* & 7 & 24 \\ \hline
Research in Algebra and Topology I \& II* & 7 & 21 \\ \hline
Research in Analysis & 2 & 7 \\ \hline
Research in Applied Mathematics I \& II & 4 & 30 \\ \hline
Research in Geometry and Linear Algebra & 2 & 13 \\ \hline
Research in Graph Theory and Combinatorics I, II, \& III & 12 & 44 \\ \hline
Research in Number Theory I \& II* & 1 & 22 \\ \hline
Teaching Introductory Mathematics* & 3 & 11 \\ \hline
Teaching Mathematics Beyond the Calculus Sequence & 3 & 8 \\ \hline
\end{tabular}
}

\newpage

\centerline{AMS Contributed Paper Sessions}\bigskip

\centerline{
\begin{tabular}{|l||c|c|}\hline
\bf Session & \bf FS & \bf TS \\ \hline\hline
Algebraic Geometry & 3 & 21 \\ \hline
Analysis and Partial Differential Equations & 4 & 13 \\ \hline
Applied Mathematics I: Mechanics, Fluids, Waves & 4 & 15 \\ \hline
Applied Mathematics II* & 4 & 11 \\ \hline
C*-Algebras and Analysis & 2 & 9 \\ \hline
Combinatorics and Number Theory & 3 & 16 \\ \hline
Combinatorics I \& II & 12 & 29 \\ \hline
Commutative Algebra and Homological Methods & 5 & 12 \\ \hline
Complex and Geometric Analysis & 3 & 11 \\ \hline
Difference Equations, Approximations, Sequences, and Special Functions & 4 & 13 \\ \hline
Differential and Integral Equations and Their Applications & 2 & 12 \\ \hline
Fractal Geometry, Complex Dynamics, and Dynamical Systems & 5 & 16 \\ \hline
Game Theory and Computing & 4 & 11 \\ \hline
Geometric Applications of Combinatorics and K-Theory & 5 & 14 \\ \hline
Geometry and General Topology** & 5 & 19 \\ \hline
Graph Theory & 4 & 12 \\ \hline
Group Theory	 & 6 & 19 \\ \hline
History of Mathematics* & 6 & 12 \\ \hline
Knots and Their Invariants & 5 & 11 \\ \hline
Knots, Topological Graphs, and Algebraic Topology & 5 & 19 \\ \hline
Lattices, Polynomials, and Linear Algebra & 3 & 13 \\ \hline
Logic and Probability & 6 & 16 \\ \hline
Mathematical Modeling and Mathematical Biology & 7 & 12 \\ \hline
Mathematics Education & 9 & 12 \\ \hline
Natural Resource Modeling and Mathematical Biology* & 7 & 15 \\ \hline
Noncommutative Algebra and Lie Theory & 2 & 13 \\ \hline
Number Theory I \& II* & 6 & 22 \\ \hline
Numerical Methods and Computing I \& II & 9 & 22 \\ \hline
Operator Theory and Banach Spaces & 6 & 16 \\ \hline
Optimization, Calculus of Variations, Nonlinear Programming & 7 & 14 \\ \hline
Partial Differential Equations & 6 & 20 \\ \hline
Probability and Stochastic Dynamical Systems & 9 & 17 \\ \hline
Statistical Modeling, Big Data, and Computing & 7 & 12 \\ \hline
Structural and Extremal Problems in Graph Theory & 3 & 6 \\ \hline
Undergraduate Research in Algebra, Combinatorics and Number Theory & 3 & 13 \\ \hline
Undergraduate Research in Analysis and Topology & 7 & 13 \\ \hline
Undergraduate Research in Applied Mathematics & 5 & 12 \\ \hline
\end{tabular}
}

\newpage

\centerline{MAA Contributed Paper Sessions}\bigskip

\centerline{
\begin{tabular}{|l||c|c||c|c|}\hline
\bf Session & \bf FS & \bf TS & \bf FO & \bf TO \\ \hline\hline
Assessing Quantitative Reasoning and Literacy & 4 & 8 & 1 & 4 \\ \hline
Assessing Student Learning: Alternative Approaches & 16 & 32 & 2 & 5 \\ \hline
Assessment of Proof Writing Throughout the Mathematics Major	 & 2 & 4 & 2 & 2 \\ \hline
Bridging the Gap: Designing an Introduction to Proofs Course & 8 & 13 & 1 & 1 \\ \hline
Data, Modeling, and Computing in the Introductory Statistics Course & 8 & 15 & 0 & 3 \\ \hline
Flipping the Classroom & 15 & 38 & 2 & 2 \\ \hline
History of Mathematical Communities	 & 4 & 10 & 2 & 2 \\ \hline
Innovative and Effective Ways to Teach Linear Algebra & 5 & 15 & 1 & 3 \\ \hline
Instructional Approaches to Increase Awareness of the Societal Value of & \multirow2*4 & \multirow2*8 & \multirow2*2 & \multirow2*2 \\
\; Mathematics &&&& \\ \hline
The Intersection of Mathematics and the Arts & 19 & 40 & 0 & 1 \\ \hline
Is Mathematics the Language of Science? & 0 & 7 & 0 & 3 \\ \hline
Mathematics Experiences in Business, Industry, and Government & 5 & 17 & 1 & 3 \\ \hline
Mathematics and Sports & 5 & 24 & 0 & 2 \\ \hline
Open Source Mathematics Textbooks & 1 & 15 & 0 & 2 \\ \hline
Programs and Approaches for Mentoring Women and Minorities in Mathematics & 7 & 8 & 2 & 2 \\ \hline
Projects, Demonstrations, and Activities that Engage Liberal Arts Mathematics & \multirow2*{13} & \multirow2*{22} & \multirow2*1 & \multirow2*1 \\
\; Students &&&& \\ \hline
Putting a Theme in a History of Mathematics Course & 4 & 9 & 0 & 2 \\ \hline
Reinventing the Calculus Sequence & 1 & 8 & 0 & 2 \\ \hline
Research on the Teaching and Learning of Undergraduate Mathematics & 6 & 19 & 1 & 3 \\ \hline
Scholarship of Teaching and Learning in Collegiate Mathematics & 9 & 15 & 2 & 5 \\ \hline
Student Activities & 9 & 20 & 2 & 2 \\ \hline
Teaching with Technology: Impact, Evaluation, and Reflection & 6 & 23 & 0 & 1 \\ \hline
Topics and Techniques for Teaching Real Analysis & 5 & 18 & 0 & 4 \\ \hline
Trends in Undergraduate Mathematical Biology Education & 4 & 12 & 0 & 1 \\ \hline
USE Math: Undergraduate Sustainability Experiences in the Introductory & \multirow2*5 & \multirow2*9 & \multirow2*2 & \multirow2*3 \\
\; Mathematics Classroom &&&& \\ \hline
Using Online Resources to Augment the Traditional Classroom & 5 & 17 & 0 & 2 \\ \hline
Wavelets in Undergraduate Education & 5 & 12 & 1 & 3 \\ \hline
We Did More with Less: Streamlining the Undergraduate Mathematics Curriculum & 3 & 5 & 1 & 2 \\ \hline
\end{tabular}
}

} 

\bigskip\bigskip

One of the aspects of the data in which we are interested is how having women among the organizers of a session is correlated with the proportion of female speakers in their sessions. J.~Bryan performed a statistical analysis on the data from this JMM that includes organizers as well as speakers, fitting a binomial regression to the data using the percentage of female organizers as the explanatory variable and the percentage of female speakers as the outcome. As is standard, the data point corresponding to each session was weighted by the number of speakers in each session.

These data points, and the fitted function, are plotted in the figure on the next page. As Bryan (personal communication) writes: ``The statistical significance of the term that capture the association between female proportion in speakers vs. organizers is extremely high. The `$z$-score' is $7.146$ and is something you can compare to a standard normal distribution. Therefore the `$p$-value'---i.e., the probability of seeing a result as or more extreme under the null of no association---is $8.91\times10^{-13}$.''

\newpage

{
\smaller\smaller

\centerline{
\includegraphics[width=5in]{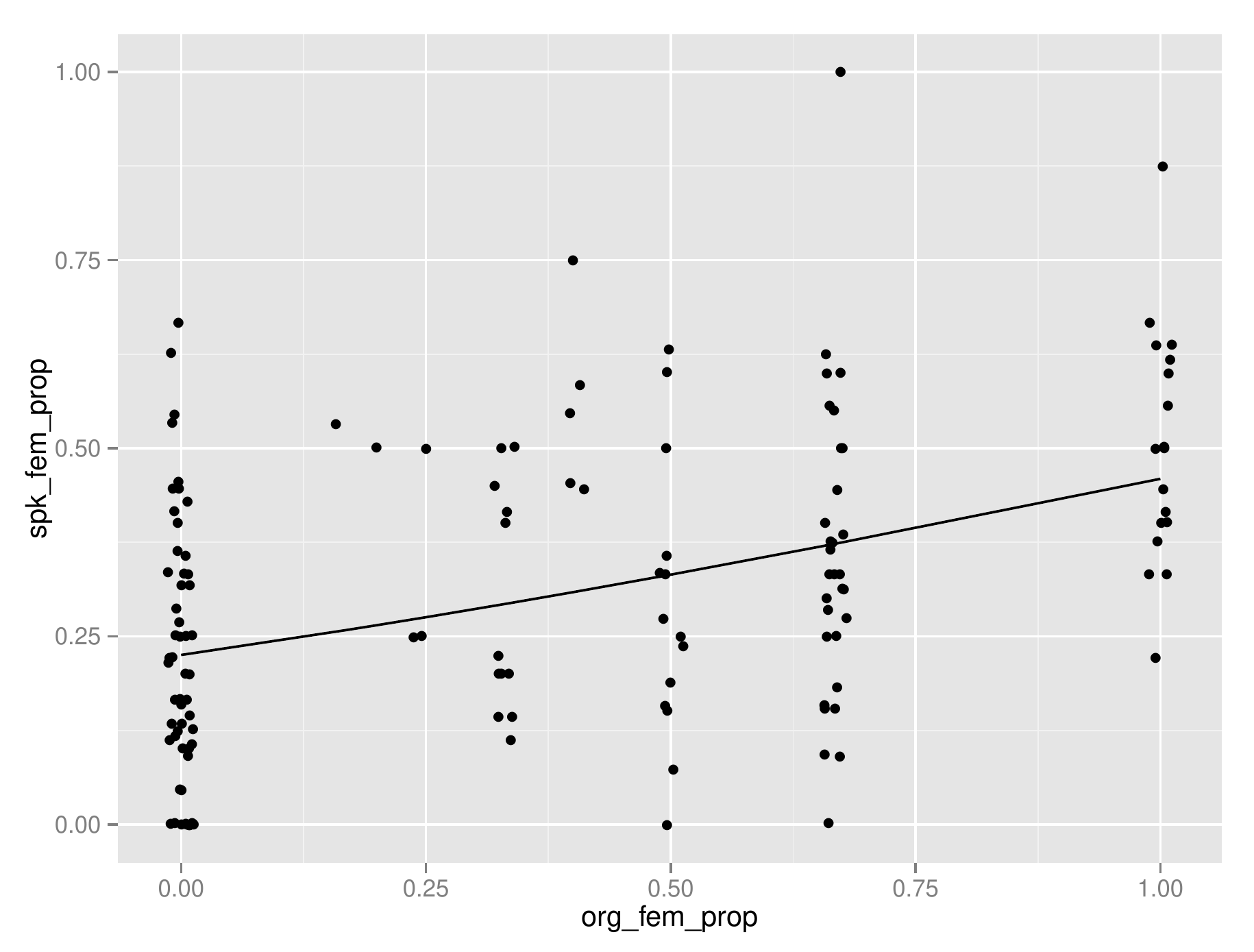}
}

\bigskip\bigskip\centerline{\bf The 2014 ICM in Seoul} \bigskip

\centerline{Invited ICM Panels}\bigskip

\centerline{
\begin{tabular}{|l||c|c||c|c|}\hline
\bf Panel & \bf FS & \bf TS & \bf FO & \bf TO \\ \hline\hline
Panel 1. Why STEM? & 1 & 3 & 1 & 1 \\ \hline
Panel 2. How should we teach better? & 0 & 3 & 1 & 1 \\ \hline
Panel 3. Mathematics is everywhere & 0 & 3 & 1 & 1 \\ \hline
Panel 4. R\&D policy (ERC+NRF)  & 0 & 4 & 0 & 1 \\ \hline
Panel 5. IMAGINARY - Panel: Math communication for the future - a Vision Slam & 1 & 4 & 0 & 1 \\ \hline
\end{tabular}
}

\bigskip\bigskip
\centerline{Panels organized by IMU}\bigskip

\centerline{
\begin{tabular}{|l||c|c||c|c|}\hline
\bf Panel & \bf FS & \bf TS & \bf FO & \bf TO \\ \hline\hline
Panel 1. Mathematical Massive Open Online Courses & 0 & 4 & 0 & 1 \\ \hline
Panel 2. Future of Publishing & 1 & 6 & 0 & 1 \\ \hline
Panel 3. World Digital Mathematics Library & 2 & 5 & 0 & 1 \\ \hline
\end{tabular}
}

\newpage

\centerline{ICM Sections}\bigskip

\centerline{
\begin{tabular}{|l||c|c|}\hline
\bf Section & \bf FS & \bf TS \\ \hline\hline
Plenary Lectures & 1 & 20 \\ \hline\hline
Algebra & 2 & 9 \\ \hline
Algebraic and Complex Geometry & 0 & 11 \\ \hline
Analysis and its Applications & 3 & 17 \\ \hline
Combinatorics  & 3 & 11 \\ \hline
Control Theory and Optimization & 1 & 6 \\ \hline
Dynamical Systems and Ordinary Differential Equations & 2 & 11 \\ \hline
Geometry* & 2 & 15 \\ \hline
History of Mathematics & 0 & 3 \\ \hline
Lie Theory and Generalizations & 1 & 11 \\ \hline
Logic and Foundations & 1 & 6 \\ \hline
Mathematical Aspects of Computer Science  & 1 & 8 \\ \hline
Mathematical Physics & 2 & 12 \\ \hline
Mathematics Education and Popularization of Mathematics (panels) & 2 & 8 \\ \hline
Mathematics Education and Popularization of Mathematics (speakers) & 1 & 3 \\ \hline
Mathematics in Science and Technology  & 5 & 10 \\ \hline
Number Theory & 1 & 15 \\ \hline
Numerical Analysis and Scientific Computing & 1 & 6 \\ \hline
Partial Differential Equations & 3 & 13 \\ \hline
Probability and Statistics & 0 & 13 \\ \hline
Topology & 0 & 10 \\ \hline
\end{tabular}
}

\bigskip\bigskip\centerline{\bf Mathematics Awards} \bigskip

\centerline{AMS Prizes}\bigskip

\centerline{
\begin{tabular}{|l||c|c|}\hline
\bf Prize & \bf FS & \bf TS \\ \hline\hline
George David Birkhoff Prize in Applied Mathematics & 1 & 16 \\ \hline
B\^ocher Memorial Prize & 0 & 33 \\ \hline
Frank Nelson Cole Prize in Algebra & 0 & 26 \\ \hline
Frank Nelson Cole Prize in Number Theory & 0 & 29 \\ \hline
Levi L. Conant Prize & 1 & 18 \\ \hline
Joseph L. Doob Prize  & 0 & 6 \\ \hline
Leonard Eisenbud Prize for Mathematics and Physics & 0 & 5 \\ \hline
Delbert Ray Fulkerson Prize** & 3 & 65 \\ \hline
E. H. Moore Research Article Prize   & 0 & 6 \\ \hline
Frank and Brennie Morgan Prize for Outstanding Research in Mathematics by an & \multirow2*3 & \multirow2*{19} \\
\; Undergraduate Student && \\ \hline
David P. Robbins Prize   & 1 & 4 \\ \hline
Ruth Lyttle Satter Prize in Mathematics   & 13 & 13 \\ \hline
Leroy P. Steele Prize for Lifetime Achievement   & 1 & 25 \\ \hline
Leroy P. Steele Prize for Mathematical Exposition  & 1 & 29 \\ \hline
Leroy P. Steele Prize for Seminal Contribution to Research  & 1 & 34 \\ \hline
Oswald Veblen Prize in Geometry   & 0 & 29 \\ \hline
Albert Leon Whiteman Memorial Prize   & 0 & 4 \\ \hline
Norbert Wiener Prize in Applied Mathematics & 0 & 13 \\ \hline
\end{tabular}
}

\newpage

\centerline{Other Prizes}\bigskip

\centerline{
\begin{tabular}{|l||c|c|}\hline
\bf Prize & \bf FS & \bf TS \\ \hline\hline
Abel Prize & 0 & 14 \\ \hline
Fields Medal & 1 & 56 \\ \hline
\end{tabular}
}

\bigskip\bigskip\centerline{\bf Editorial Boards of Mathematics Journals} \bigskip

} 

These ten journals were selected by consulting multiple lists of the ``top'' mathematics journals. The list of journals was finalized before looking at any of their editorial boards, to avoid biasing the sample. The editorial boards (as of October 2014) were then located on the internet, and the genders of their members determined via internet searches.

\smaller\smaller

\bigskip\bigskip

\centerline{
\begin{tabular}{|l||c|c|}\hline
\bf Journal & \bf FS & \bf TS \\ \hline\hline
Acta Mathematica & 0 & 9 \\ \hline
Annals of Mathematics & 0 & 12 \\ \hline
Annales Scientifiques de l'\'Ecole Normale Sup\'erieure & 2 & 9 \\ \hline
Communications on Pure and Applied Mathematics & 0 & 13 \\ \hline
Duke Mathematical Journal & 1 & 22 \\ \hline
Inventiones Mathematicae & 0 & 10 \\ \hline
Journal f\"ur die Reine und Angewandte Mathematik & 0 & 5 \\ \hline
Journal of Differential Geometry & 0 & 15 \\ \hline
Journal of the American Mathematical Society & 4 & 26 \\ \hline
Proceedings of the London Mathematical Society & 4 & 45 \\ \hline
\end{tabular}
}

\end{document}